# FullSWOF:
# A SOFTWARE FOR OVERLAND FLOW SIMULATION
*FullSWOF : un logiciel pour la simulation du ruissellement*

**Olivier DELESTRE**

Laboratoire J.A. Dieudonné & Polytech Nice – Sophia
Postal address[1]: UMR CNRS 7351 UNSA, 06108 Nice Cedex 02, France
e-mail : delestre@unice.fr

**Stéphane CORDIER**

Laboratoire MAPMO – Fédération Denis Poisson
Postal address: UMR CNRS 7349, Université d'Orléans, B.P. 6759, 45067 Orléans Cedex 2, France
e-mail : stephane.cordier@math.cnrs.fr

**Frédéric DARBOUX**

Institut National de la Recherche Agronomique (INRA)
Postal address: UR 0272 Science du sol, Centre de recherche d'Orléans, CS 40001, 45075 Orléans Cedex 2, France
e-mail : Frederic.Darboux@orleans.inra.fr

**Mingxuan DU**

Polytech Nice – Sophia
Postal address: 930 Route des Colles – BP 145, 06903 Sophia Antipolis Cedex, France
e-mail : mingxuanism@gmail.com

**François JAMES**

Laboratoire MAPMO – Fédération Denis Poisson
Postal address: UMR CNRS 7349, Université d'Orléans, B.P. 6759, 45067 Orléans Cedex 2, France
e-mail : Francois.James@math.cnrs.fr

**Christian LAGUERRE**

Laboratoire MAPMO – Fédération Denis Poisson
Postal address: UMR CNRS 7349, Université d'Orléans, B.P. 6759, 45067 Orléans Cedex 2, France
e-mail : christian.laguerre@math.cnrs.fr

**Carine LUCAS**

Laboratoire MAPMO – Fédération Denis Poisson
Postal address: UMR CNRS 7349, Université d'Orléans, B.P. 6759, 45067 Orléans Cedex 2, France
e-mail : carine.lucas@univ-orleans.fr

**Olivier PLANCHON**

Institut de Recherche pour le Développement (IRD)
Postal address: UMR LISAH, 2 Place Viala, 34060 Montpellier Cedex 1, France
e-mail : olivier.planchon@gmail.com

*Le ruissellement sur les terres agricoles peut avoir des effets indésirables tels que l'érosion des sols, les inondations et le transport de polluants. Afin de mieux comprendre ce phénomène et d'en limiter les conséquences, nous avons développé un code à l'aide de méthodes numériques récentes : FullSWOF (Full Shallow Water equations for Overland Flow) un code orienté objet écrit en C++. Il est libre et peut être téléchargé à partir de http://www.univ-orleans.fr/mapmo/soft/FullSWOF/. Le modèle est basé sur le système de Saint-Venant. Les difficultés numériques viennent des nombreuses transitions sec/mouillé et de la topographie très variable rencontrée sur le terrain. Le code intègre le ruissellement, les précipitations, l'infiltration (Green-Ampt modifié), la friction (les lois de Darcy-Weisbach et de Manning).*
*Nous présentons d'abord la méthode numérique pour la résolution des équations en eaux peu profondes integrée dans FullSWOF_2D (la version en deux dimensions). Cette méthode est basée sur le schéma de reconstruction hydrostatique, couplée à un traitement semi-implite du terme de friction. FullSWOF_2D a déjà été validé à l'aide des*

---

[1] Corresponding author

*solutions analytiques de la bibliothèque SWASHES. FullSWOF_2D est exécuté sur des données de terrain acquises sur une parcelle située à Thiès (Sénégal). Les résultats de la simulation sont comparés avec les données mesurées. Ce banc d'essai expérimental permet de démontrer les capacités de FullSWOF à simuler l'écoulement de surface. FullSWOF pourrait également être utilisé pour d'autres problèmes environnementaux, tels que les inondations fluviales et les ruptures de barrage.*

*Overland flow on agricultural fields may have some undesirable effects such as soil erosion, flood and pollutant transport. To better understand this phenomenon and limit its consequences, we developed a code using state-of-the-art numerical methods: FullSWOF (Full Shallow Water equations for Overland Flow), an object oriented code written in C++. It has been made open-source and can be downloaded from http://www.univ-orleans.fr/mapmo/soft/FullSWOF/. The model is based on the classical system of Shallow Water (SW) (or Saint-Venant system). Numerical difficulties come from the numerous dry/wet transitions and the highly-variable topography encountered inside a field. The code includes runon and rainfall inputs, infiltration (modified Green-Ampt equation), friction (Darcy-Weisbach and Manning formulas).*
*First we present the numerical method for the resolution of the Shallow Water equations integrated in FullSWOF_2D (the two-dimension version). This method is based on hydrostatic reconstruction scheme, coupled with a semi-implicit friction term treatment. FullSWOF_2D has been previously validated using analytical solutions from the SWASHES library (Shallow Water Analytic Solutions for Hydraulic and Environmental Studies). FullSWOF_2D is run on a real topography measured on a runoff plot located in Thies (Senegal). Simulation results are compared with measured data. This experimental benchmark demonstrate the capabilities of FullSWOF to simulate adequately overland flow. FullSWOF could also be used for other environmental issues, such as river floods and dam-breaks.*

## Key words

## I   INTRODUCTION

Rain on agricultural fields can yield to the occurrence of overland flow. At field scale, overland flow may have some undesirable effects such as soil erosion and pollutant transport. Downstream the watersheds, roads and houses may be damaged. Some control measures can be taken such as using grass strips. We have to know how the water is moving in order to put these developments. In overland flow prediction, several methods are used from black box models to physically based models. Two physical models are often used to model overland flow: Kinematic (KW) and Diffusive Wave (DW) equations [14,15]. But following [24,10,11], we choose to use the Shallow Water (de Saint Venant [21]) physical model. Indeed KW and DW models may give poor results in terms of water heights and velocities in case of mixed subcritical and supercritical flow. Inspite of computational time SW is mandatory. MacCormack scheme is widely used to solve SW equations [24,10,11]. But it neither guarantees the positivity of water depths at the wet/dry transitions, nor preserves steady states (not well-balanced [13]) as noticed in [20]. In industrial codes (ISIS, Canoe, HEC-RAS, MIKE11...), SW equations are often solved under non-conservative form [15] with either Preissmann scheme or Abbott-Ionescu scheme. Thus transcritical flows and hydraulic jump are not solved properly. In order to cope with all these problems, we choose to use the hydrostatic reconstruction [1,2]. This positive preserving well-balanced finite volume scheme is integrated in FullSWOF_2D. In what follows, we present the physical model, then the numerical methods. In the end, FullSWOF_2D is applied on a real event measured in Thies by IRD [23].

## II   THE MODEL

### II.1   The shallow water equations

As in [10,11], we consider the 2D Shallow Water equations (SW2D) which write (see **Figure 1 a**)

$$\partial_t h + \partial_x(hu) + \partial_y(hv) = R - I$$
$$\partial_t(hu) + \partial_x(hu^2 + gh^2/2) + \partial_y(huv) = gh(S_{0x} - S_{fx}) \quad , \tag{1}$$
$$\partial_t(hv) + \partial_x(huv) + \partial_y(hv^2 + gh^2/2) = gh(S_{0y} - S_{fy})$$

where the unknowns are the velocities $u(x,y,t)$ and $v(x,y,t)$ [m/s] and the water height $h(x,y,t)$ [m]. The subscript $x$ (respectively $y$) stands for the x-direction (resp. the y-direction): $S_{0x} = -\partial_x z(x,y)$ and

$S_{0y}=-\partial_y z(x,y)$ are the ground slopes, $S_{fx}$ and $S_{fy}$ are friction terms. $R(x,y,t)$ [m/s] is the rainfall intensity and $I(x,y,t)$ [m/s] is the infiltration rate. As in [10], we use the Darcy-Weisbach friction law which writes:

$$S_{fx}=f\frac{u\sqrt{u^2+v^2}}{8gh}, \quad S_{fy}=f\frac{v\sqrt{u^2+v^2}}{8gh} \quad . \tag{2}$$

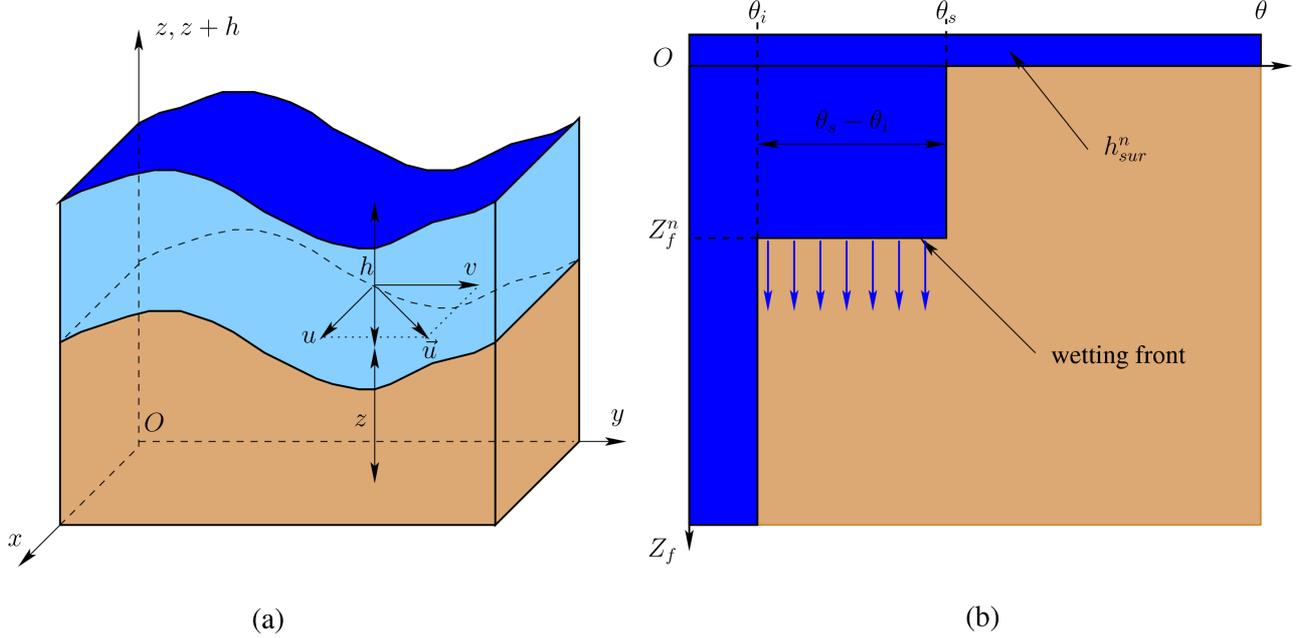

**Figure 1**: Illustration of variables of (a) Shallow Water equations (SW2D) and (b) Green-Ampt infiltration model.

**II.2    The Green-Ampt infiltration model**

Infiltration is computed at each cell using a modified version of the Green-Ampt model [12,18,10]. With this model, the movement of water in soil is assumed to be in the form of an advancing wetting front (located at $Z_f^n$ [m]) that separates a zone still at the initial soil moisture $\theta_s$ (see **Figure 1 b**). At the moment $t=t_n$, the infiltration capacity $I_C^n$ [m/s] is calculated thanks to

$$I_C^n = K_s\left(1+\frac{h_f - h_{sur}^n}{Z_f^n}\right) \quad \text{where} \quad Z_f^n = \frac{V_{inf}^n}{\theta_s - \theta_i}, \tag{3}$$

where $h_f$ is the wetting front capillary pressure head, $K_s$ the hydraulic conductivity at saturation, $h_{sur}^n$ the water height and $V_{inf}^n$ the infiltrated water volume. Thus we have the infiltration rate

$$I^n = \frac{min(h_{sur}^n, \Delta t . I_C^n)}{\Delta t} \quad , \tag{4}$$

and the infiltrated volume

$$V_{inf}^{n+1} = V_{inf}^n + \Delta t . I^n \quad , \tag{5}$$

where $\Delta t$ is the time step. In the case of a two-layer soil, we consider a modification of this model (see [5,10]).

**III THE NUMERICAL METHOD**

The scheme will be presented in one dimension (SW1D)

$$\begin{aligned}\partial_t h + \partial_x(hu) &= R - I \\ \partial_t(hu) + \partial_x(hu^2 + gh^2/2) &= gh(S_{0x} - S_{fx})\end{aligned} , \tag{6}$$

it is integrated in FullSWOF_1D. The one dimensional version of FullSWOF has been developed for teaching, academical research and numerical development purposes. The extension of the 1D numerical method to SW2D on structured grid is straightforward and is integrated in FullSWOF_2D (Full Shallow-Water equations for Overland Flow). FullSWOF_2D is an object-oriented C++ code (free software and GPL-compatible license CeCILL-V2, source code available at http://www.univ-orleans.fr/mapmo/soft/FullSWOF/) developed in the framework of the project ANR METHODE (for details about the code see [4,5]). FullSWOF is structured in order to ease the implementation of new numerical methods and physical models. In what follows, we note the discharge $q=hu$ [m²/s] and the vector of conservative variables $U=\begin{pmatrix} h & hu \end{pmatrix}^t$.

### III.1 Convective step

A finite volume discretization of SW1D writes

$$U_i^* = U_i^n - \frac{\Delta t}{\Delta x}\left[F_{i+1/2L}^n - F_{i-1/2R}^n - Fc_i^n\right], \qquad (7)$$

with $\Delta x$ the space step and

$$F_{i+1/2L}^n = F_{i+1/2}^n + S_{i+1/2L}^n$$
$$F_{i-1/2R}^n = F_{i-1/2}^n + S_{i-1/2R}^n, \qquad (8)$$

the left (respectively right) modification of the numerical flux $Fl$ for the homogeneous problem (see section III.3)

$$F_{i+1/2}^n = Fl(U_{i+1/2L}^n, U_{i+1/2R}^n). \qquad (9)$$

The values $U_{i+1/2L}$ and $U_{i+1/2R}$ are obtained thanks to two consecutive reconstructions. Firstly a MUSCL reconstruction [2,5] is performed on $u$, $h$ and $h+z$ in order to get a second order scheme in space. This gives us the reconstructed values $(U_-, z_-)$ and $(U_+, z_+)$. Secondly we apply the hydrostatic reconstruction [1,2] on the water height which allows us to get a positive preserving well-balanced scheme (in the sense it preserves at least steady state at rest)

$$\begin{aligned} h_{i+1/2L} &= max(h_{i+1/2-} + z_{i+1/2-} - max(z_{i+1/2-}, z_{i+1/2+}), 0) \\ U_{i+1/2L} &= (h_{i+1/2L}, h_{i+1/2L} u_{i+1/2-}) \\ h_{i+1/2R} &= max(h_{i+1/2+} + z_{i+1/2+} - max(z_{i+1/2-}, z_{i+1/2+}), 0) \\ U_{i+1/2R} &= (h_{i+1/2R}, h_{i+1/2R} u_{i+1/2+}) \end{aligned} \qquad (10)$$

We introduce

$$S_{i+1/2L}^n = \begin{pmatrix} 0 \\ \frac{g}{2}(h_{i+1/2-}^2 - h_{i+1/2L}^2) \end{pmatrix}, \quad S_{i-1/2R}^n = \begin{pmatrix} 0 \\ \frac{g}{2}(h_{i-1/2+}^2 - h_{i-1/2R}^2) \end{pmatrix} \qquad (11)$$

and a centered source term is added to preserve consistency and well-balancing (see [1,2])

$$Sc_i = \begin{pmatrix} 0 \\ -g\frac{h_{i-1/2+} + h_{i+1/2-}}{2}(z_{i+1/2-} - z_{i-1/2+}) \end{pmatrix}. \qquad (12)$$

We have to insist on the positivity and the robustness of this method. But, let us mention that such method requires some condition on the mesh size compared with mean water height (see [6]). The rain and the infiltration are treated explicitly (for details see [5]).

### III.2 Friction treatment

In this step, the friction term is taken into account with the following system

$$\partial_t U = \begin{pmatrix} 0 \\ -ghS_f \end{pmatrix}. \qquad (13)$$

This system is solved thanks to a semi-implicit method (as in [3,11])

$$\begin{aligned} h^{n+1} &= h^* \\ q^{n+1} &= \frac{q^*}{1 + \Delta t \frac{f}{8}\frac{|q^n|}{h^n h^{n+1}}} \end{aligned} \qquad (14)$$

where $h^*$, $q^*$ and $u^*$ are the variables from the convective step. This method allows to preserve stability (under a classical CFL condition) and steady states at test. Finally, these two steps are combined in a second order TVD Runge Kutta method which is the Heun's predictor-corrector method. It writes

$$U^* = U^n + \Delta t\, \Phi(U^n)$$
$$U^{**} = U^* + \Delta t\, \Phi(U^*)$$
$$U^{n+1} = \frac{U^n + U^{**}}{2}, \qquad (15)$$

where $\Phi$ is the right part of (7).

### III.3 Numerical flux

We use the HLL flux [2,5] which writes

$$Fl(U_L, U_R) = \begin{cases} F(U_L) & \text{if } 0 < c_1 \\ \dfrac{c_2 F(U_L) - c_1 F(U_R)}{c_2 - c_1} + \dfrac{c_1 c_2}{c_2 - c_1}(U_R - U_L) & \text{if } c_1 < 0 < c_2 \\ F(U_R) & \text{if } c_2 < 0 \end{cases}, \qquad (16)$$

with two parameters $c_1 < c_2$ given by

$$c_1 = \min_{U = U_L, U_R}\left(\min_{j \in \{1,2\}} \lambda_j(U)\right), \quad c_2 = \max_{U = U_L, U_R}\left(\max_{j \in \{1,2\}} \lambda_j(U)\right), \qquad (17)$$

where $\lambda_1(U) = u - \sqrt{gh}$ and $\lambda_2(U) = u + \sqrt{gh}$ are the eigenvalues of SW1D. In practice, we use a CFL condition $n_{CFL} = 0.5$ at second order and $n_{CFL} = 1$ at first order, with

$$\Delta t \leq n_{CFL} \frac{\Delta x}{\max_{i \in [\![1,J]\!]} (|u_i| + \sqrt{gh_i})}, \qquad (18)$$

where $J$ is the number of space cells. At second order, variables $(h_i, u_i)$ in (18) are replaced by the reconstructed values $(h_{i+1/2-}, u_{i+1/2-})$ and $(h_{i+1/2+}, u_{i+1/2+})$ (detailed in next section).

### III.4 MUSCL- reconstruction

We define the MUSCL reconstruction of a scalar function $s \in \mathbb{R}$ [22] by

$$s_{i-1/2+} = s_i - \frac{\Delta x}{2} Ds_i, \quad s_{i+1/2-} = s_i + \frac{\Delta x}{2} Ds_i, \qquad (19)$$

with the operator

$$Ds_i = \text{minmod}\left(\frac{s_i - s_{i-1}}{\Delta x}, \frac{s_{i+1} - s_i}{\Delta x}\right), \qquad (20)$$

and the *minmod* limiter

$$\text{minmod}(x,y) = \begin{cases} \min(x,y) & \text{if } x,y \geq 0 \\ \max(x,y) & \text{if } x,y \leq 0 \\ 0 & \text{else} \end{cases}. \qquad (21)$$

As mentioned previously, the MUSCL reconstruction is performed on $u$, $h$ and $h+z$ then we deduce the reconstruction of $z$. In order to keep the discharge conservation, the reconstruction of the velocity has to be modified as what follows

$$u_{i-1/2+} = u_i - \frac{h_{i+1/2-}}{h_i} \frac{\Delta x}{2} Du_i, \quad u_{i+1/2-} = u_i + \frac{h_{i-1/2+}}{h_i} \frac{\Delta x}{2} Du_i. \qquad (22)$$

## IV VALIDATION

FullSWOF_2D has already been validated on analytical solutions integrated in SWASHES (Shallow Water Analytic Solutions for Hydraulic and Environmental Studies) a free library of analytical solutions written in object-oriented ISO C++ [7,8]. The purpose of this section is to confront FullSWOF_2D to a real system: the plot of Thies, Senegal [23]. The results presented below are the results preliminary to a more detailed study. They simply aim at illustrating the ability of FullSWOF_2D to simulate a dynamic runoff on real data without any preprocessing.

## IV.1 Experimental data

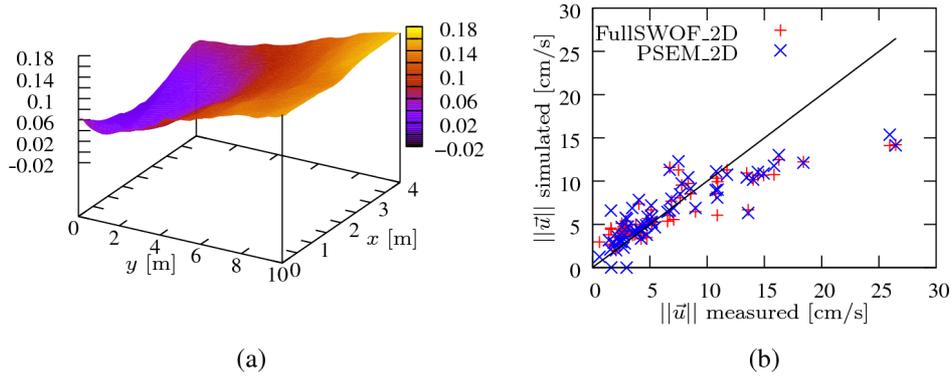

(a)                                                          (b)

**Figure 2**: (a) The topography of the experimental plot and (b) comparison between measured and simulated velocities.

Thies (Senegal) is the plot of an experimental system instrumented by IRD [23] in the project PNRH RIDES. It consisted of a 2-hour artificial rainfall on a sandy-soil plot of 4-by-10 meters (**Figure 2 a**). The plot had the classical configuration of Wooding's open book, with 1% slope along (Ox) and (Oy) axes. Several experiments have been carried out to test the Salt Velocity Gauge [16] (a new measurement technique) and to study the dynamics of runoff and erosion. The set of measures considered consists of flow measurement at the outlet and the flow velocity measured at 63 locations in the plot. These data have been used to compare different computer codes: NCF and MAHLERAN (based on KW), RillGrow (based on DW) and PSEM_2D (based on SW2D with MacCormack scheme). The dataset is freely available at http://www.umr-lisah.fr/Thies_2004/. In the following, we have used the set of parameters obtained with PSEM_2D and we have observed the sensitivity to infiltration parameters.

## IV.2 Numerical results

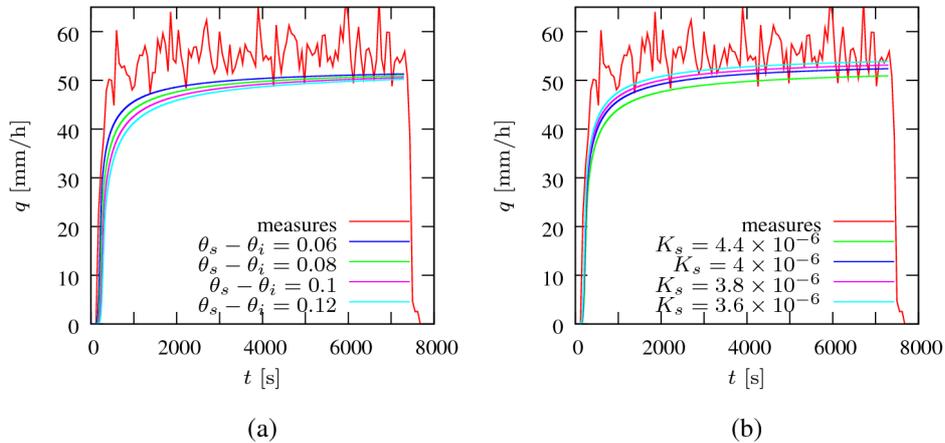

(a)                                                          (b)

**Figure 3**: Comparison between simulated hydrographs and measures (a) for several values of $\theta_s - \theta_i$ and (b) $K_s$.

For all the simulations, we have used $\Delta x = \Delta y = 0.1\,m$ as space steps and $f = 0.26$ as Darcy Weisbach friction coefficient. The average rain intensity was $70\,mm/h$ during two hours. Infiltration parameters used with PSEM_2D [23] were $h_f = 0.06\,m$, $\theta_s - \theta_i = 0.12$ and $K_s = 4.4\text{e}-6\,m/s$. We tried this set of parameters. We notice that the simulated velocities with FullSWOF_2D are closed to those obtained with PSEM_2D (**Figure 2 b**). Small flow-velocities are well caught, excepted for the bigger values. In these cases, values are underestimated both by FullSWOF_2D and PSEM_2D. This means that the friction law is not adapted to do this simulation. Moreover we notice that the simulated hydrograph is under the measured one. Thus we have tried other values for $\theta_s - \theta_i$ and $K_s$. We kept $h_f$ and $\theta_s - \theta_i$ (respectively $K_s$) values and we changed $K_s$ (resp. $\theta_s - \theta_i$). We notice that decreasing $\theta_s - \theta_i$ improves mainly the beginning of the hydrograph (**Figure 3 a**), while decreasing $K_s$ improves the entire hydrogram (**Figure 3 b**). We notice that the saturated hydraulic conductivity $K_s$ is the most influent parameter, as observed in [19] a study on the impact of this parameter on the surface runoff. Moreover, we have performed the simulation of this event on the integrated hydrological modelling system: MIKE SHE by DHI. This software is widely used for overland flow simulation, thus it might be used as a reference. In this software, the overland flow is calculated using a finite difference method using the diffusive wave approximation coupled with the one-dimensional Richards' equation [17] for the infiltration. Thus it is an approximation of the physical model included in FullSWOF_2D. We notice that FullSWOF_2D gives water heights very closed to those obtained with MIKE SHE (**Figure 4**).

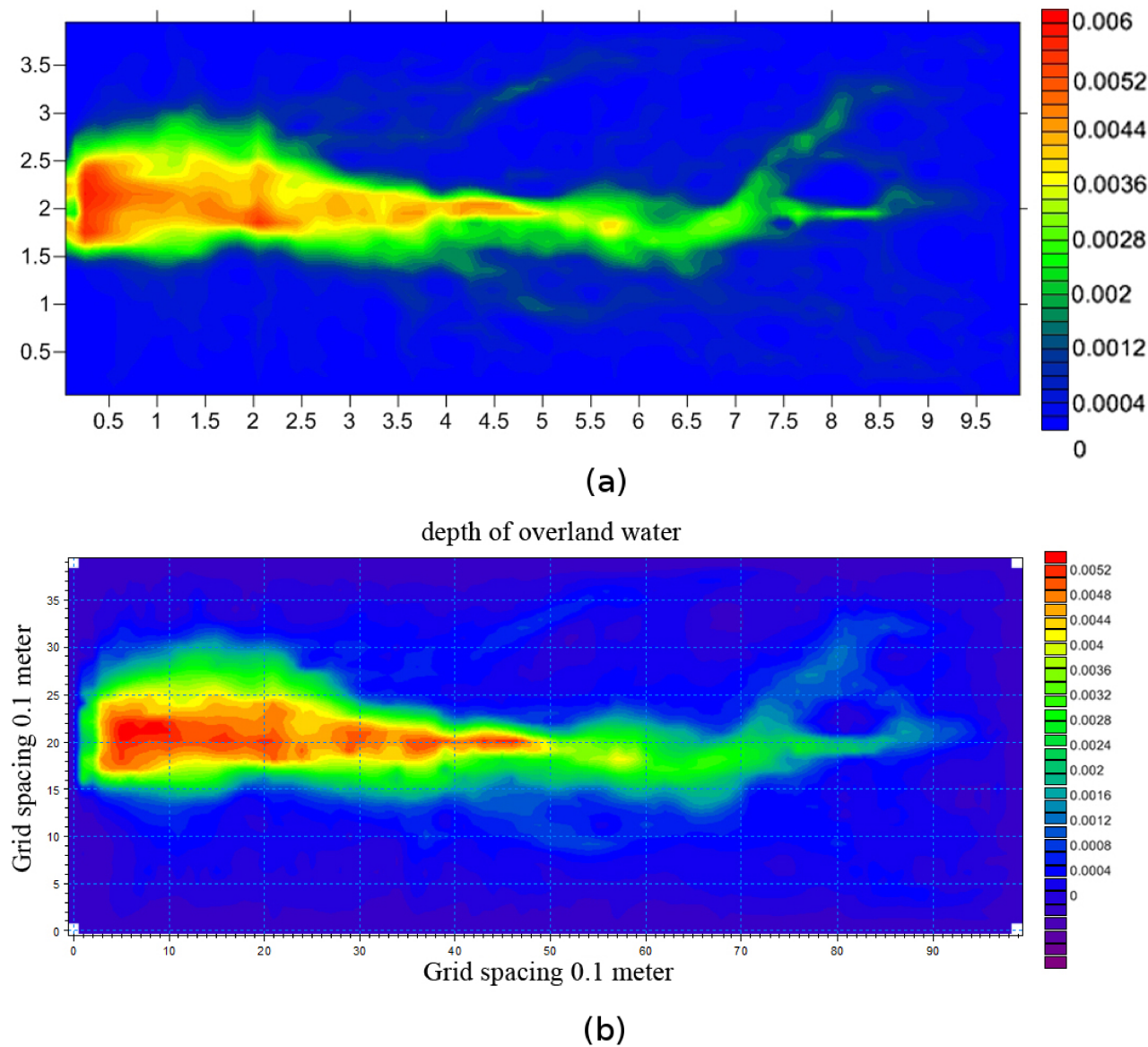

**Figure 4**: Comparison between the water heights obtained with (a) FullSWOF_2D and (b) MIKE SHE.

This work is still in progress, we will have to do more complete comparison between FullSWOF, measured data and other softwares.

## V CONCLUSIONS

FullSWOF is a freely-available object oriented code designed to simulate overland flow on agricultural fields. This code is based on the Shallow Water system. The numerical difficulties such as dry/wet transitions and steady states are dealt with using up-to-date numerical methods. FullSWOF gave good results on real data. Nevertheless, the physical model has to be improved to catch more accurately the measured velocities.

## VI ACKNOWLEGMENTS AND THANKS

This work was partially supported by ANR grant "METHODE"#ANR-07-BLAN-0232.

# VII REFERENCES AND CITATIONS

[1] Audusse, E., Bouchut, F., Bristeau, M.-O., Klein, R., & Perthame, B. (2004). A fast and stable well-balanced scheme with hydrostatic reconstruction for shallow water flows. *J. Sci. Comp.*, **25(6)**, 2050-2065.

[2] Bouchut, F. (2004). *Nonlinear stability of finite volume methods for hyperbolic conservation laws, and well-balanced schemes for sources*. Frontiers in Mathematics, Birkhauser.

[3] Bristeau, M.-O., & Coussin, B. (2001). *Boundary conditions for the shallow water equations solved by kinetic schemes*. Inria report RR-4282.

[4] Delestre, O. (2008). *Ecriture d'un code C++ pour la simulation en hydrologie*. Master thesis Université d'Orléans, available on http://dumas.ccsd.cnrs.fr/dumas-00446163/fr/.

[5] Delestre, O. (2010). *Simulation du ruissellement d'eau de pluie sur des surfaces agricoles*. PhD thesis Université d'Orléans, available on http://tel.archives-ouvertes.fr/INSMI/tel-00531377/fr.

[6] Delestre, O., Cordier, S., Darboux, F. & James, F. (2012). A limitation of the hydrostatic reconstruction technique for Shallow Water equations, *C. R. Acad. Sci. Paris, Ser. I*, http://dx.doi.org/10.1016/j.crma.2012.08.004.

[7] Delestre, O., Lucas, C., Ksinant, P.-A., Darboux, F., Laguerre, C., Vo, T.N.T., James, F., & Cordier, S. (submitted). SWASHES: a compilation of Shallow Water Analytic Solutions for Hydraulic and Environmental Studies, available on http://hal.archives-ouvertes.fr/hal-00628246.

[8] Delestre, O., Lucas, C., Ksinant, P.-A., Darboux, F., Laguerre, C., James, F., & Cordier, S. (submitted). SWASHES: a library for benchmarking in hydraulics/SWASHES : Une bibliothèque de bancs d'essai en hydraulique. To be published in *Proceedings of SimHydro 2012,* Polytech'Nice Sophia, Sophia-Antipolis (France), 2012. Advances in Hydroinformatics – SimHydro 2012, Springer.

[9] El Bouajaji, M. (2007). *Modélisation des écoulements à surface libre : étude du ruissellement des eaux de pluie*. Master thesis Université Louis Pasteur, available on http://dumas.ccsd.cnrs.fr/dumas-00459336/fr/.

[10] Esteves, M., Faucher, X., Galle, S., & Vauclin, M. (2000). Overland flow and infiltration modelling for small plots during unsteady rain: numerical results versus observed values. *J. Hydrol.*, **228**, 265-282.

[11] Fiedler, F. R., & Ramirez, J. A. (2000). A numerical method for simulating discontinuous shallow flow over an infiltrating surface. *Int. J. Numer. Methods Fluids*, **32**, 219-240.

[12] Green, W. H., & Ampt, G. A. (1911). Studies on soil physics. *The Journal of Agricultural Science*, **4**, 1-24.

[13] Greenberg, J. M., & LeRoux, A.-Y. (1996). A well-balanced scheme for the numerical processing of source terms in hyperbolic equation. *SIAM Journal on Numerical Analysis*, **33**, 1-16.

[14] Moussa, R., & Bocquillon, C. (2000). Approximation of the Saint-Venant equations for flood routing with overbank flow. *SIAM Hydrology and Earth System Sciences*, **4**, 251-261.

[15] Novak, P., Guinot, V., Jeffrey, A., & Reeve, D.E. (2010). *Hydraulic modelling – an Introduction*. Spoon Press.

[16] Planchon, O., Silvera, N., Gimenez, R., Favis-Mortlock, D., Wainwright, J., Le Bissonnais, Y., & Govers, G. (2005). An automated salt-tracing gauge for flow-velocity measurement. E*arth Surface Processes and Landforms*, **30**, 833-844.

[17] Richards, L.A. (1931). Capillary conduction of liquids through porous mediums. *Physics*, **1**, 318-333.

[18] Rousseau, M. (2008). *Modélisation des écoulements à surface libre : étude du ruissellement des eaux de pluie*. Master thesis Université de Nantes, available on http://dumas.ccsd.cnrs.fr/dumas-00494243/fr/.

[19] Rousseau, M., Cerdan, O., Ern, A., Le Maître, O., & Sochala, P. (2012). A study of overland flow with uncertain infiltration using stochastic tools. *Advances in Water Resources*, **38(0)**, 1-12.

[20] Rousseau, M., Cerdan, O., Delestre, O., Dupros, F., James, F., & Cordier, S. (submitted). Overland flow modelling with Shallow Water Equation using a well balanced numerical scheme: Adding efficiency or just more complexity ?, available on http://hal.archives-ouvertes.fr/hal-00664535.

[21] de Saint-Venant, A. J.-C. (1871). Théorie du mouvement non-permanent des eaux, avec application aux crues des rivières et à l'introduction des marées dans leur lit. *Comptes Rendus de l'Académie des Sciences*, **73**, 147-154.

[22] van Leer, B. (1979). Towards the ultimate conservative difference scheme. V. A seond-order sequel to Godunov's method. *Journal of Computational Physics*, **32(1)**, 101-136.

[23] Tatard, L., Planchon, O., Wainwright, J., Nord, G., Favis-Mortlock, D., Silvera, N., Ribolzi, O., Esteves, M. & Huang, C.-H. (2008). Measurement and modelling of high-resolution flow-velocity data under simulated rainfall on a low-slope sandy soil. *Journal of Hydrology*, **348**, 1-12.

[24] Zhang, W., & Cundy, T.W. (1989). Modeling of two-dimensional overland flow. *Water Resources Research*, **25**, 2019-2035.